 \documentclass[preprint,review,12pt]{elsarticle}
\usepackage[margin=1.5 in]{geometry}

\usepackage{amssymb}

\makeatletter
\def\ps@pprintTitle{%
 \let\@oddhead\@empty
 \let\@evenhead\@empty
 \def\@oddfoot{\centerline{\thepage}}%
 \let\@evenfoot\@oddfoot}
\makeatother

\usepackage{xcolor}

\usepackage{graphics} 
\usepackage{epsfig} 
\usepackage{amsmath} 
\usepackage{amsthm}
\usepackage{amssymb}  
\usepackage{subfigure}

\begin{document}

\begin{frontmatter}


\title{Comment on "Event-Triggered Stabilization of Linear Time-Delay Systems via Halanay-Type Inequality"\tnoteref{label0}}
\tnotetext[label0]{This work was supported by the Natural Sciences and Engineering Research Council of Canada (NSERC) under Discovery Grant RGPIN-2022-03144. K. Zhang further acknowledges support from the Verner Family Faculty Fellowship in Smith Engineering at Queen's University.}
\author[QU1]{Luke Jonker}\ead{23jhf4@queensu.ca}

\author[QU1,QU2]{Kexue Zhang}\ead{kexue.zhang@queensu.ca}

\address[QU1]{Department of Mathematics and Statistics, Queen's University, Canada}
\address[QU2]{Department of Electrical and Computer Engineering, Queen's University, Canada}

\begin{abstract}
This comment revisits Lemma~1 in~\cite{KZ2023}, which plays a central role in the event-triggered stabilization analysis developed therein. We identify technical gaps in the proof of the lemma and provide a corrected argument. In particular, careful treatment of the exponentially decaying term shows that its decay rate must be retained in the resulting convergence estimate. The statement of the original lemma, with the exponential decay rate determined by the minimum of the characteristic decay rate and the decay rate of this term, remains valid.
\end{abstract}


\end{frontmatter}

We consider Lemma~1 in~\cite{KZ2023}. The following argument provides a corrected proof while preserving the statement of the original lemma. The technical gaps in the original proof and how the argument below addresses them are discussed after the proof. Throughout this comment, we use the same notation as in~\cite{KZ2023}.

\medskip

\noindent\textbf{Lemma~1}~\cite{KZ2023}.
Let $a,b,\alpha,\beta,r$ be positive constants satisfying
\[
a>b+\alpha.
\]
Let $v:[-r,\infty)\rightarrow\mathbb{R}^{+}$ be a continuous function satisfying
\begin{equation}
D^{+}v(t)
\leq
-av(t)
+b\sup_{-r\leq s\leq0}\{v(t+s)\}
+\alpha\|v_0\|_{r}e^{-\beta t},
\label{eq:main}
\end{equation}
for $t\in\mathbb{R}^{+}$, where
\[
\|v_0\|_{r}
=
\sup_{-r\leq s\leq0}\{v(s)\}.
\]
Then
\begin{equation}
v(t)\leq \|v_0\|_{r}e^{-\eta t},
\qquad t\in\mathbb{R}^{+},
\label{eq:conclusion}
\end{equation}
where
\[
\eta=\min\{\lambda,\beta\},
\]
and $\lambda>0$ is the unique solution of
\begin{equation}
a=be^{\lambda r}+\lambda+\alpha.
\label{eq:lambda}
\end{equation}

\medskip

\noindent\textit{Proof.}
Since $a>b+\alpha$, there exists a sufficiently small $\varepsilon>0$ such that
\[
a-\varepsilon>b+\alpha.
\]
Let $\lambda_{\varepsilon}>0$ be the unique solution of
\begin{equation}
a-\varepsilon
=
be^{\lambda_{\varepsilon}r}
+\lambda_{\varepsilon}
+\alpha,
\label{eq:lambda-eps}
\end{equation}
and define
\begin{equation}
\eta_{\varepsilon}
:=
\min\{\lambda_{\varepsilon},\beta\}.
\label{eq:eta-eps}
\end{equation}
For an arbitrary $\xi>0$, set
\[
B:=\|v_0\|_{r}+\xi
\]
and define
\begin{equation}
w(t):=v(t)-Be^{-\eta_{\varepsilon}t},
\qquad t\geq-r.
\label{eq:w}
\end{equation}

For $-r\leq t\leq0$, we have
\[
v(t)
\leq
\|v_0\|_{r}
<
B
\leq
Be^{-\eta_{\varepsilon}t},
\]
where the last inequality follows from $t\leq0$ and
$\eta_{\varepsilon}>0$. Hence,
\begin{equation}
w(t)<0,
\qquad -r\leq t\leq0.
\label{eq:initial-w}
\end{equation}
In particular,
\[
w(0)=v(0)-B
\leq
\|v_0\|_{r}-(\|v_0\|_{r}+\xi)
=-\xi<0.
\]

We next show that $w(t)<0$ for all $t>0$. Suppose, to the contrary, that there exists some $t>0$ such that $w(t)\geq0$. Define
\[
t^{*}
:=
\inf\{t\geq0:w(t)\geq0\}.
\]
Since $w(0)<0$ and $w$ is continuous, it follows that
$t^{*}>0$. Moreover, continuity of $w$ gives
\begin{equation}
w(t^{*})=0,
\qquad
w(t)<0,\quad 0\leq t<t^{*}.
\label{eq:first-contact}
\end{equation}
Thus,
\begin{equation}
v(t^{*})=Be^{-\eta_{\varepsilon}t^{*}},
\label{eq:contact}
\end{equation}
and
\[
v(t)<Be^{-\eta_{\varepsilon}t},
\qquad
0\leq t<t^{*}.
\]
Together with the estimate on $[-r,0]$, this yields
\[
v(t)\leq Be^{-\eta_{\varepsilon}t},
\qquad
-r\leq t\leq t^{*}.
\]

Consequently, for $0\leq t\leq t^{*}$,
\begin{align}
\sup_{-r\leq s\leq0}\{v(t+s)\}
&\leq
\sup_{-r\leq s\leq0}
\left\{
Be^{-\eta_{\varepsilon}(t+s)}
\right\}
\nonumber\\
&=
e^{\eta_{\varepsilon}r}
Be^{-\eta_{\varepsilon}t}.
\label{eq:sup}
\end{align}

We next justify a lower estimate for $v$ near $t=t^*$. Define
\[
q(t)
:=
\frac{v(t)}
{Be^{-\eta_{\varepsilon}t}}.
\]
Since $v$ is continuous and the denominator is strictly positive,
$q$ is continuous. By \eqref{eq:contact},
\[
q(t^{*})=1.
\]
Since $a>\varepsilon$, we have
\[
0<1-\frac{\varepsilon}{a}<1.
\]
Therefore, by continuity of $q$, there exists a sufficiently small
$\delta>0$, with $t^{*}-\delta\geq0$, such that
\[
q(t)\geq1-\frac{\varepsilon}{a},
\qquad
t\in(t^{*}-\delta,t^{*}).
\]
Equivalently,
\begin{equation}
v(t)
\geq
\left(1-\frac{\varepsilon}{a}\right)
Be^{-\eta_{\varepsilon}t},
\qquad
t\in(t^{*}-\delta,t^{*}).
\label{eq:lower}
\end{equation}

Using \eqref{eq:main}, \eqref{eq:sup}, and \eqref{eq:lower}, for
$t\in(t^{*}-\delta,t^{*})$ we obtain
\begin{align}
D^{+}w(t)
={}&
D^{+}v(t)
+\eta_{\varepsilon}Be^{-\eta_{\varepsilon}t}
\nonumber\\
\leq{}&
-av(t)
+b\sup_{-r\leq s\leq0}\{v(t+s)\}
\nonumber\\
&\quad
+\alpha\|v_0\|_{r}e^{-\beta t}
+\eta_{\varepsilon}Be^{-\eta_{\varepsilon}t}
\nonumber\\
\leq{}&
\left[
-a+\varepsilon
+be^{\eta_{\varepsilon}r}
+\eta_{\varepsilon}
\right]
Be^{-\eta_{\varepsilon}t}
\nonumber\\
&\quad
+\alpha\|v_0\|_{r}e^{-\beta t}.
\label{eq:key1}
\end{align}

At this point, it is essential to retain the exponentially decaying term. Since
\[
\eta_{\varepsilon}
=
\min\{\lambda_{\varepsilon},\beta\}
\leq\beta,
\]
we have, for $t\geq0$,
\[
e^{-\beta t}
\leq
e^{-\eta_{\varepsilon}t}.
\]
Moreover,
\[
\|v_0\|_{r}\leq B.
\]
It therefore follows that
\begin{equation}
\alpha\|v_0\|_{r}e^{-\beta t}
\leq
\alpha Be^{-\eta_{\varepsilon}t}.
\label{eq:forcing}
\end{equation}
Combining \eqref{eq:key1} and \eqref{eq:forcing} gives
\begin{align}
D^{+}w(t)
\leq{}&
\left[
-a+\varepsilon
+be^{\eta_{\varepsilon}r}
+\eta_{\varepsilon}
+\alpha
\right]
Be^{-\eta_{\varepsilon}t}.
\label{eq:key2}
\end{align}

Since $\eta_{\varepsilon}\leq\lambda_{\varepsilon}$ and the function
\[
x\mapsto be^{xr}+x
\]
is increasing for $x\geq0$, we have
\[
be^{\eta_{\varepsilon}r}+\eta_{\varepsilon}
\leq
be^{\lambda_{\varepsilon}r}
+\lambda_{\varepsilon}.
\]
Hence, by~\eqref{eq:key2} and~\eqref{eq:lambda-eps},
\begin{align}
D^{+}w(t)
&\leq
\left[
-a+\varepsilon
+be^{\lambda_{\varepsilon}r}
+\lambda_{\varepsilon}
+\alpha
\right]
Be^{-\eta_{\varepsilon}t}
\nonumber\\
&=0
\label{eq:dini}
\end{align}
for all $t\in(t^{*}-\delta,t^{*})$.

Since $w$ is continuous and its upper Dini derivative is nonpositive throughout $(t^{*}-\delta,t^{*})$, the standard monotonicity property of Dini derivatives implies that $w$ is nonincreasing on this interval (see, e.g.,~\cite{VL-XZL:1993}). Consequently,
\[
w(t^{*})
\leq
w(t^{*}-\delta)
<0,
\]
which contradicts $w(t^{*})=0$.

Therefore,
\[
w(t)<0,
\qquad t\geq0,
\]
and hence
\begin{equation}
v(t)
<
(\|v_0\|_{r}+\xi)e^{-\eta_{\varepsilon}t},
\qquad t\geq0.
\label{eq:eps-bound}
\end{equation}

It remains to pass to the limit as $\varepsilon\to0$. Define
\[
F(x):=be^{xr}+x+\alpha,
\qquad x\geq0.
\]
The function $F$ is continuous and strictly increasing. Equations
\eqref{eq:lambda} and \eqref{eq:lambda-eps} give
\[
F(\lambda)=a,
\qquad
F(\lambda_{\varepsilon})=a-\varepsilon.
\]
Thus, by continuity of the inverse of the strictly increasing function
$F$,
\[
\lambda_{\varepsilon}\to\lambda
\qquad\text{as}\qquad
\varepsilon\to0.
\]
Consequently,
\[
\eta_{\varepsilon}
=
\min\{\lambda_{\varepsilon},\beta\}
\longrightarrow
\min\{\lambda,\beta\}
=\eta.
\]
Letting first $\varepsilon\to0$ and then $\xi\to0$ in
\eqref{eq:eps-bound} yields
\[
v(t)\leq\|v_0\|_{r}e^{-\eta t},
\qquad t\geq0,
\]
which proves the result.
\hfill$\square$

\medskip

We now clarify the technical issues in the proof of Lemma~1 in~\cite{KZ2023}. First, there is a gap in the estimate
 \[
 \sup_{-r\leq s\leq 0}\left\{v(t^*+s)\right\}
 < \sup_{-r\leq s\leq 0}\left\{\varepsilon \left(\|v_0\|_r +\xi \right) e^{-\eta (t^*+s)}\right\}
 \]
because equality is attained when $s=0$ according to the definition of $t^*$. Moreover, the contradiction argument in~\cite{KZ2023} is not valid because $D^{+}w(t^*)<0$ does not imply $w$ is strictly decreasing in the interval $(t^*-\delta,t^*+\delta)$ for small enough $\delta > 0$. A simple counterexample is $w(t)=-|t|$, and $D^{+}w(0)=-1<0$ but $w$ is strictly increasing on $(-\infty,0)$.

The argument presented here avoids both issues by establishing the estimate~\eqref{eq:sup} directly and deriving $D^+w(t)\le 0$ for $t\in (t^*-\delta,t^*)$. Since the upper Dini derivative is nonpositive throughout the interval $(t^*-\delta,t^*)$, the standard monotonicity property of Dini derivatives implies that $w$ is nonincreasing on this interval. This yields the desired contradiction and completes the proof. We emphasize that these issues concern only the proof of Lemma~1 in~\cite{KZ2023}; the statement of the lemma remains valid.


\begin{thebibliography}{99}

\bibitem{KZ2023}
K. Zhang,
``Event-triggered stabilization of linear time-delay systems via Halanay-type inequality".
\emph{IEEE Control Syst. Lett.}, vol. 7, pp. 3205--3210, 2023.

\bibitem{VL-XZL:1993}
V. Lakshmikantham and X. Z. Liu, 
\emph{Stability Analysis in Terms of Two
Measures}. 
Singapore: World Sci., 1993.
\end{thebibliography}
\end{document}